\documentclass[11pt]{article}

\usepackage{cite}
\usepackage{latexsym,color,amsmath,amsthm,amssymb,amscd,amsfonts}
\usepackage{fancyhdr,graphicx}

\def\be{\begin{equation}}
\def\ee{\end{equation}}
\def\bea{\begin{eqnarray}}
\def\eea{\end{eqnarray}}
\def\bt{\begin{theorem}}
\def\et{\end{theorem}}
\def\bl{\begin{lemma}}
\def\el{\end{lemma}}
\def\br{\begin{remark}}
\def\er{\end{remark}}
\def\bc{\begin{corollary}}
\def\ec{\end{corollary}}
\def\bd{\begin{definition}}
\def\ed{\end{definition}}

\def\b{\beta}

\def\l{\lambda}

\def\vP{\vec{P}}
\def\vQ{\vec{Q}}
\def\vf{\vec{f}}

\def\bP{\mathbf{P}}
\def\bQ{\mathbf{Q}}

\def\bbR{\mathbb{R}}

\def\b1{B_{1}^}

\def\dfpq{D_{\vec{\mathbf{f}}} (\vec{\bP}, \vec{\bQ})}
\def\ba{\begin{array}}
\def\ea{\end{array}}
\def\ben{\begin{enumerate}}
\def\een{\end{enumerate}}
\newtheorem{theorem}{Theorem}[section]
\newtheorem{lemma}{Lemma}[section]
\newtheorem{remark}{Remark}[section]
\newtheorem{proposition}{Proposition}[section]
\newtheorem{corollary}{Corollary}[section]

\newtheorem{definition}{Definition}[section]
\def\dfpqi{D_{\vf}(\vP, \vQ; i)}
\def\bff{\mathbf{f}}

\begin{document}

\title{Mixed $f$-divergence for multiple pairs of measures\footnote{Keywords: Alexandrov-Fenchel  inequality, $\bf f$-dissimilarity, $f$-divergence, isoperimetric inequality.
Mathematics Subject Classification Number: 28, 52, 60.}}  
\author{Elisabeth M. Werner \thanks{Partially supported by an  NSF grant} and Deping Ye \thanks{Supported by an NSERC grant and a
start-up grant from Memorial University of Newfoundland} }
\date{}
\maketitle
 
\begin{abstract}
In this paper,  the concept of the classical $f$-divergence for a pair of measures is extended to the mixed $f$-divergence for multiple pairs of measures. The mixed $f$-divergence provides a way to measure the difference between multiple pairs of (probability) measures. Properties for the mixed $f$-divergence are established, such as permutation invariance and symmetry in distributions. An
Alexandrov-Fenchel type inequality and an isoperimetric inequality for the 
mixed $f$-divergence are proved. 
\end{abstract}

\section{Introduction} 

In applications such as pattern matching, image analysis, statistical learning, and information theory, one often needs to compare two (probability) measures and needs to know whether they are similar to each other. Hence, finding the ``right" quantity to measure the difference between two (probability) measures $P$ and $Q$ is central. Traditionally, people use the classical $L_p$ distances between $P$ and $Q$, such as the variational distance and the $L_2$ distance. However,  the family of $f$-divergences is often more suitable to fulfill the goal than the classical $L_p$ distance of measures. 
\par
The $f$-divergence $D_f(P, Q)$ of two probability measures $P$ and $Q$ was first introduced in \cite{Csiszar} and independently in \cite{AliSilvery1966, Morimoto1963} and was defined by
 \begin{equation}
 D_f(P, Q)=\int_Xf\left(\frac{p}{q}\right) q\,d\mu. \label{f:divergence}
\end{equation}
Here,  $p$ and $q$ are  density functions of  $P$ and $Q$ with respect to  a measure $\mu$ on $X$.
The idea behind the $f$-divergence is to replace, for instance, the function $f(t)=|t-1|$ in the variational distance by a general convex function $f$. Hence the $f$-divergence includes various widely used divergences as special cases, such as, the variational distance, the  Kullback-Leibler  divergence \cite{KullbackLeibler1951}, the Bhattacharyya distance \cite{Bhattacharyya1946} and many more. Consequently,  the $f$-divergence receives considerable attention  
 not only in the information theory (e.g.,  \cite{BarronGyorfiMeulen, CoverThomas, HarremoesTopsoe, LieseVajda2006, OsterrVajda}) but also in many other areas. We only mention convex geometry. Within the last few years,  amazing connections have been discovered between notions and  concepts  from convex geometry and information theory, e.g., \cite{Gardner2002, GLYZ, JenkinsonWerner, LutwakYangZhang2004/1, LutwakYangZhang2005, PaourisWerner2011}, 
leading to a totally new point of view and introducing  a whole new set of tools in the area of convex geometry. 
In particular, it was observed in \cite{Werner2012/1} that one of the most important affine invariant notions,
the $L_p$-affine surface area for convex bodies, e.g.,  \cite{Ludwig2010, Ludwig-Reitzner, Ludwig-Reitzner1999, Lu1, SW2004},  is R\'enyi entropy from information theory and statistics.   R\'enyi entropies are special cases of $f$-divergences and consequently those were then introduced for convex bodies and their corresponding entropy inequalities have been established in \cite{Werner2012b}. We also refer to,  for instance \cite{Basseville2010},  for more references related to the $f$-divergence.

\par
Extension of the $f$-divergence from two (probability) measures to multiple (probability) measures is fundamental in many applications,  such as statistical hypothesis test and classification, and much  research has been devoted to that, for instance in   \cite{  Menendez, MoralesPardo1998, Zografos1998}. Such extensions include, e.g., the Matusita's affinity \cite{Matusita1967, Matusita1971}, the Toussaint's affinity \cite{Toussaint1974},  the information radius \cite{Sibson1969} and the average divergence \cite{Sgarro1981}. 
\par
 The $\bf f$-dissimilarity $D_{\bff}(P_1, \cdots, P_l)$ for (probability) measures $P_1, \cdots, P_l$,  introduced in \cite{GyorfiNemetz1975, GyorfiNemetz1978} for a convex function  $\bff : \bbR^l\rightarrow \bbR$, is   
 a natural generalization of the $f$-divergence. 
 It is defined as
 \begin{equation*} 
D_\bff(P_1, \cdots, P_l)=\int_X \bff(p_1, \cdots, p_l)\,d\mu, \end{equation*}  where the $p_i$'s are density functions of the $P_i$'s that are absolutely continuous with respect to $\mu$. 
For a convex function $f$, the function $\bff(x,y)=y f(\frac{x}{y})$ is also convex on $x, y>0$, and $D_{\bff}(P, Q)$ is equal to the classical $f$-divergence defined in formula (\ref{f:divergence}). Note that the Matusita's affinity is related to  $$\bff(x_1,\cdots, x_l)=-\prod_{i=1}^lx_i^{1/l},$$ and the Toussaint's affinity is related to  $\bff(x_1,\cdots, x_l)=-\prod_{i=1}^lx_i^{a_i}$, where $ a_i\geq 0$ and  such that $\sum_{i=1}^la_i=1.$ 
\par
Here, we introduce  special   $\bf f$-dissimilarities, namely the mixed $f$-divergence and the $i$-th mixed $f$-divergence, which can be viewed as vector forms of the usual $f$-divergence. We establish some basic properties of these quantities,  such as permutation invariance and symmetry in distributions. We prove an isoperimetric type inequality and  an Alexandrov-Fenchel type inequality for the 
mixed $f$-divergence. Alexandrov-Fenchel  inequality is a fundamental inequality in convex geometry  and  many important inequalities such as  the Brunn-Minkowski inequality and  Minkowski's first inequality
follow from it (see, e.g., \cite{Gardner2002, Sch}).
\par
The paper is organized as follows.  In Section 2 we establish some basic properties of the mixed $f$-divergence, such as permutation invariance and symmetry in distributions. 
 In Section 3 we prove the 
general Alexandrov-Fenchel inequality and  isoperimetric  inequality for the mixed $f$-divergence.
Section 4 is dedicated to the $i$-th mixed $f$-divergence and its related isoperimetric  type inequalities.

\section{The Mixed $f$-Divergence}

Throughout this paper, let $(X, \mu)$ be a finite measure space. For $1 \leq i
\leq n$, let  $ P_i=p_i  \mu$ and  $ Q_i=q_i \mu$ be probability
measures on $X$ that are absolutely continuous with respect to
the measure $\mu$. {Moreover, we assume that for all $i=1,  \cdots, n$,  $p_i$ and $q_i$ are nonzero $\mu$-a.e.} We use $\vec{\bP}$ and  $\vec{\bQ}$ to denote the vectors of probability measures, or, in short,   probability vectors,  $$\vec{\bP}=(P_1, P_2, \cdots, P_n), \ \ \ \vec{\bQ}=(Q_1, Q_2, \cdots, Q_n).$$ We  use $\vec{p}$ and $\vec{q}$ to denote the vectors of density functions, or density vectors,  for $\vec{\bP}$ and $\vec{\bQ}$ respectively, $$\frac{\,d\vec{\bP}}{\,d\mu}= \vec{p}=(p_1, p_2, \cdots, p_n),\ \ \ \ \frac{\,d\vec{\bQ}}{\,d\mu} =\vec{q}=(q_1, q_2, \cdots, q_n).$$ 
We make the convention that $0 \cdot \infty =0$. 
 
Denote by $\mathbb{R}^+=\{x\in \mathbb{R}: x\geq 0\}$.  Let $f: (0, \infty) \rightarrow  \mathbb{R}^+$ be a non-negative convex or concave function. 
The $*$-adjoint function $f^*:(0, \infty) \rightarrow  \mathbb{R}^+$
of $f$  is defined by
\begin{equation*}
f^*(t) = t f (1/t).
\end{equation*}
It is obvious that $(f^*)^*=f$ and that $f^*$ is again convex,  respectively concave, if
$f$ is convex,  respectively concave.
 
\par
 Let $f_i: (0, \infty) \rightarrow
\mathbb{R}^+$, $ 1 \leq i \leq n$, be either convex or concave functions. Denote by $\vec{\mathbf{f}} =(f_1, f_2, \cdots, f_n)$ the vector of functions. We write $$\vec{\mathbf{f}}^*=(f_1^*, f_2^*, \cdots, f_n^*)$$ to be the $*$-adjoint vector for $\vec{\mathbf{f}}$. 
\vskip 2mm 
Now  we introduce {\em the mixed $f$-divergence} for $(\vec{\mathbf{f}}, \vec{\bP}, \vec{\bQ})$  as follows. 
\vskip 2mm 
\begin{definition}\label{mixedf} Let  $(X, \mu)$ be a measure space. Let $\vec{\bP}$ and $\vec{\bQ}$ be two probability  vectors on $X$ with density vectors $\vec{p}$ and $\vec{q}$ respectively.  
The mixed $f$-divergence $\dfpq$ for $(\vec{\mathbf{f}}, \vec{\bP}, \vec{\bQ})$ 
is defined by
\begin{equation}\label{mixed1} 
\dfpq= \int_{X} \prod_{i=1}^n
\left[f_i\left(\frac{p_{i}}{q_{i}}\right) q_{i}\right]^\frac{1}{n}
d \mu.
\end{equation}
\end{definition}
\par
\noindent 
Similarly, we  define the mixed $f$-divergence for $(\vec{\mathbf{f}}, \vec{\bQ}, \vec{\bP})$ by
\begin{equation}\label{mixed2}
D_{\vec{\mathbf{f}}}(\vec{\bQ}, \vec{\bP})=  \int_{X} \prod_{i=1}^n
\left[f_i\left(\frac{q_{i}}{p_{i}}\right) p_{i}
\right]^\frac{1}{n} d\mu.
\end{equation}
A  special case is when all distributions $P_i$ and $Q_i$ are identical and equal to a probability distribution $P$. 
In this case, \begin{eqnarray*} \dfpq=D_{(f_1, f_2, \cdots, f_n)}\big((P, P, \cdots, P),(P, P, \cdots, P)\big)= \prod_{i=1}^n\left[f_i(1)\right]^{\frac{1}{n}}.\end{eqnarray*} 
\vskip 2mm
Let $\pi\in S_n$ denote a permutation on $\{1, 2, \cdots, n\}$ and denote $$\pi(\vec{p})=(p_{\pi(1)}, p_{\pi(2)},\cdots, p_{\pi(n)}).$$
One immediate result from  Definition \ref{mixedf} is the following permutation invariance for $\dfpq$. 
\par
\begin{proposition}
[\bf Permutation invariance] Let the vectors $\vec{\mathbf{f}}, \vec{\bP}, \vec{\bQ}$ be as above, and let $\pi\in S(n)$ be a permutation on $\{1, 2, \cdots, n\}$. Then  $$\dfpq=D_{\pi(\vec{\mathbf{f}})}(\pi(\vec{\bP}),\pi(\vec{\bQ})).$$ 
\end{proposition}
\vskip 2mm  
When all $(f_i, P_i, Q_i)$ are equal to
$(f, P, Q)$, the mixed $f$-divergence is equal to the classical $f$-divergence, denoted by $D_f(P, Q)$, which
takes the  form
\begin{eqnarray*} D_f(P, Q)&=&D_{(f, f, \cdots, f)}\big((P, P, \cdots, P), (Q, Q, \cdots, Q)\big)=\int_{X}
f\left(\frac{p}{q}\right) q d \mu.\end{eqnarray*}
\vskip 2mm
As $f^*(t)=tf(1/t)$, one  easily obtains a fundamental property for the classical $f$-divergence $D_f(P, Q)$, namely,   $$D_f(P, Q)=D_{f^*}(Q, P),$$  for all $(f, P, Q)$. Similar results hold true for  the mixed $f$-divergence. We show this now.
\par 
Let $0\leq k\leq n$. We write $D_{\vec{\mathbf{f}}, k}(\vec{\bP}, \vec{\bQ})$ for 
\begin{eqnarray*}
D_{\vec{\mathbf{f}}, k}(\vec{\bP}, \vec{\bQ})   =     \int_{X} \prod_{i=1}^k
\left[f_i  \left(  \frac{p_{i}}{q_{i}}  \right)   q_{i}\right]^\frac{1}{n}    \times     \prod_{i=k+1}^n  
\left[f_i^*  \left(  \frac{q_{i}}{p_{i}}  \right)   p_{i}\right]^\frac{1}{n}    
d \mu.\end{eqnarray*}
 Clearly, $D_{\vec{\mathbf{f}}, n}(\vec{\bP}, \vec{\bQ})=\dfpq$ and $D_{\vec{\mathbf{f}}, 0}(\vec{\bP}, \vec{\bQ})=D_{\vec{\mathbf{f}}^*}(\vec{\bQ}, \vec{\bP})$, where $$\vec{\mathbf{f}}^*=(f_1^*, f_2^*, \cdots, f_n^*).$$
Then we have the following result for  changing order of distributions.  
\begin{proposition}[\bf Principle for changing order of distributions] \label{principle}  Let $\vec{\mathbf{f}}, \vec{\bP}, \vec{\bQ}$ be as above.  Then, for any $0\leq k\leq n$, one has $$\dfpq=D_{\vec{\mathbf{f}}, k}(\vec{\bP}, \vec{\bQ}).$$
In particular, $$\dfpq=D_{\vec{\mathbf{f}}^*}(\vec{\bQ}, \vec{\bP}).$$ 
\end{proposition} 
\vskip 2mm \noindent {\bf Proof.} Let $0\leq k\leq n$. Then, \begin{eqnarray*}
\dfpq        &=&               \int_{X} \prod_{i=1}^k
\left[f_i\left(  \frac{p_{i}}{q_{i}}  \right)   q_{i}\right]^\frac{1}{n}      \times       \prod_{i=k+1}^n    
\left[f_i\left(  \frac{p_{i}}{q_{i}}  \right)   q_{i}\right]^\frac{1}{n}
    d \mu\\&=&               \int_{X} \prod_{i=1}^k
\left[f_i\left(  \frac{p_{i}}{q_{i}}  \right)   q_{i}\right]^\frac{1}{n}      \times       \prod_{i=k+1}^n    
\left[f_i^*\left(  \frac{q_{i}}{p_{i}}  \right) p_{i}\right]^\frac{1}{n}        
d \mu \\ &=&       D_{\vec{\mathbf{f}}, k}(\vec{\bP}, \vec{\bQ}),
\end{eqnarray*} where the second equality follows from $f_i\left(\frac{p_{i}}{q_{i}}\right) q_{i}=f_i^*\left(\frac{q_{i}}{p_{i}}\right) p_{i}$.  
\vskip 2mm
A direct consequence of Proposition \ref{principle} is the following symmetry principle for the mixed $f$-divergence. 

\begin{proposition} [\bf Symmetry in distributions] Let $\vec{\mathbf{f}}, \vec{\bP}, \vec{\bQ}$ be as above. Then,  $\dfpq+D_{\vec{\mathbf{f}}^*}(\vec{\bP}, \vec{\bQ})$ is symmetric in $\vec{\bP}$ and $\vec{\bQ}$, namely, $$\dfpq+D_{\vec{\mathbf{f}}^*}(\vec{\bP}, \vec{\bQ})=  D_{\vec{\mathbf{f}} }(\vec{\bQ}, \vec{\bP})+D_{\vec{\mathbf{f}}^*}(\vec{\bQ}, \vec{\bP}).$$
\end{proposition} 
\par
\noindent {\bf Remark.} Proposition \ref{principle} says that $\dfpq$ remains the same if one replaces any triple $(f_i, P_i, Q_i)$ by $(f_i^*, Q_i, P_i)$. It is also easy to see that, for all $0\leq k, l\leq n$, one has  $$\dfpq= D_{\vec{\mathbf{f}}, k}(\vec{\bP}, \vec{\bQ})=D_{\vec{\mathbf{f}}^*,l}(\vec{\bQ}, \vec{\bP})=D_{\vec{\mathbf{f}}^*}(\vec{\bQ}, \vec{\bP}).$$ Hence, for all $0\leq k, l\leq n$,  $$D_{\vec{\mathbf{f}}, k}(\vec{\bP}, \vec{\bQ})+D_{\vec{\mathbf{f}}^*,l}(\vec{\bP}, \vec{\bQ})=\dfpq+D_{\vec{\mathbf{f}}^*}(\vec{\bP}, \vec{\bQ})$$ is symmetric in $\vec{\bP}$ and $\vec{\bQ}$. 

\vskip 2mm
 Hereafter, we only consider the  mixed $f$-divergence $\dfpq$ defined in  formula (\ref{mixed1}).  Properties for the mixed $f$-divergence $D_{\vec{\mathbf{f}}}(\vec{\bQ}, \vec{\bP})$  defined in (\ref{mixed2})  follow along the same lines. 
 \par
Now we list  some important mixed $f$-divergences.

\vskip 2mm
\noindent
{\bf Examples.}
\par
\noindent (i) The total variation is a widely used $f$-divergence to measure the difference between two probability measures $P$ and $Q$ on $(X,\mu)$. It is related to function $f(t)=|t-1|$. Similarly,  the {\em mixed total variation} is defined by $$D_{TV}(\vec{\bP}, \vec{\bQ})=\int_X \prod_{i=1}^n |p_i-	q_i|^{\frac{1}{n}}\,d\mu.$$ It  measures the difference between two probability vectors $\vec{\bP}$ and $\vec{\bQ}$. 
\par \noindent
(ii) For $a \in \mathbb{R}$, we denote by $a_+=\max\{a,0\}.$ The {\em mixed relative entropy} or {\em mixed Kullback Leibler divergence} of  $\vec{\bP}$ and $\vec{\bQ}$ is defined by
\begin{eqnarray*}
D_{KL}\big(\vec{\bP}, \vec{\bQ})   =   D_{(f_+, \cdots, f_+)}\big(\vec{\bP}, \vec{\bQ})   =      \int_{X} \prod_{i=1}^n
\bigg[ p_i \ln\bigg(\frac{q_{i}}{p_{i}}  \bigg)   
\bigg]_+^\frac{1}{n}     d\mu,
\end{eqnarray*} where $f(t) = t \ln t$. 
When $ P_i=P=p \mu$ and $ Q_i= Q =q\mu$ for all $i=1, 2, \cdots, n$, we get the following  (modified) {\em relative entropy} or {\em  Kullback Leibler divergence} $$ D_{KL}\big(P || Q\big)= \int_{X} p \left[\ln\left(\frac{q}{p}\right)\right]_+ 
d\mu. $$
\par
\noindent
(iii) For the (convex and/or 
concave) functions $f_{\alpha_i}(t) = t^{\alpha_i}$, $\alpha_i \in \mathbb{R}$ for $1 \leq i \leq n$,  the {\em mixed Hellinger integrals} is defined by 
\begin{eqnarray*}
D_{  (f_{\alpha_1}, f_{\alpha_2}, \cdots, f_{\alpha_n}) }\big(\vec{\bP}, \vec{\bQ})=\int_{X} \prod_{i=1}^n
 \left[ p_i  ^\frac{\alpha_i}{n} q_i ^\frac{1-\alpha_i}{n}\right] 
 d\mu.
\end{eqnarray*}
In particular,  $$D_{  (t^{\alpha }, t^{\alpha }, \cdots, t^{\alpha }) }\big(\vec{\bP}, \vec{\bQ}) = \int_{X} \prod_{i=1}^n
 p_i  ^\frac{\alpha}{n} q_i ^\frac{1-\alpha}{n} 
 d\mu.$$ Those integrals are related to the Toussaint's affinity \cite{Toussaint1974}, and can be used to define the {\em mixed $\alpha$-R\'enyi divergence}
 \begin{eqnarray*}
D_{\alpha}\big(\{P_i || Q_i\}_{i=1}^n\big) &=& \frac{1}{\alpha -1}  \ln \left( \int_{X} \prod_{i=1}^n
 p_i  ^\frac{\alpha}{n} q_i ^\frac{1-\alpha}{n} 
 d\mu \right)\\ &=&\frac{1}{\alpha -1}  \ln \big[D_{  (t^{\alpha }, t^{\alpha }, \cdots, t^{\alpha }) }\big(\vec{\bP}, \vec{\bQ}) \big].
\end{eqnarray*}
The case $\alpha_i=\frac{1}{2}$, for all $i=1, 2, \cdots, n$,  gives the  {\em mixed Bhattacharyya  coefficient} or {\em mixed Bhattacharyya distance} of $(\vec{\bP}, \vec{\bQ})$,  
 \begin{eqnarray*}
D_{  \big(\sqrt{t}, \sqrt{t}, \cdots, \sqrt{t}\big) }\big(\vec{\bP}, \vec{\bQ}) = \int_{X} \prod_{i=1}^n
 p_i  ^\frac{1}{2n} q_i ^\frac{1}{2n} 
 d\mu.
\end{eqnarray*} This integral is related to the Matusita's affinity  \cite{Matusita1967, Matusita1971}. For more information on the corresponding $f$-divergences we refer to e.g. \cite{LieseVajda2006}. 

\par
\noindent (iv)  In view of existing connections  between  information theory and
convex geometry (e.g.,  \cite{PaourisWerner2011, Werner2012/1, Werner2012b}), we define the mixed $f$-divergences for convex bodies (convex and compact subsets in $\mathbb R^n$ with nonempty interiors) $K_i$ with positive curvature functions $f_{K_i}$, $1\leq i\leq n$,  is via the measures 
\begin{equation*} 
\,d P_{K_i}=  \frac{1}{h_{K_i}^n}\,d\sigma \ \ \ \mathrm{and}  \ \ \ \,d Q_{K_i} =  {f_{K_i} h_{K_i}}\,d\sigma,  \ \ \ \ \ 1 \leq i \leq n.
\end{equation*}  Here, $\sigma$ is the spherical measure of the  unit sphere $S^{n-1}$,  $h_K(u) =\max_{x\in K}\langle x, u \rangle$ 
is the support function of $K$, and $f_K(u)$ is the curvature function of $K$ at $u\in S^{n-1}$, the reciprocal of the Gauss curvature at $x$ on the boundary of $K$ with unit outer normal $u$.   If $f_i: (0, \infty) \rightarrow \mathbb{R}^+$, $ 1 \leq i \leq n$, are convex and/or concave functions, then  \begin{eqnarray*}
&& D_{\vec{\mathbf{f}}}\big((P_{K_1},  \dots, P_{K_n}), (Q_{K_1}, \dots, Q_{K_n})\big) =  \int_{S^{n-  1}}    \prod_{i=1}^n 
\bigg [f_i  \bigg(  \frac{1}{ f_{K_i}h_{K_i}^{n+1}}     \bigg)    {f_{K_i} h_{K_i}} \bigg]^\frac{1}{n} \, d\sigma, \end{eqnarray*}  are the general mixed affine surface areas introduced in 
\cite{Ye2012}.  We refer to \cite{Sch} for more details on convex bodies.

\section{Inequalities}
 
The classical Alexandrov-Fenchel
inequality for  mixed volumes of convex bodies is a fundamental result in (convex) geometry. 
A general version of this inequality  for  {\em mixed volumes} of convex bodies 
can be found in  \cite{ Ale1937, Bus1958,
Sch}.  Alexandrov-Fenchel type inequalities for 
 (mixed) affine surface areas can be found in
\cite{Lut1987, Lu1, WernerYe2010, Ye2012}. Now we
prove  an inequality for the 
mixed $f$-divergence for measures, which we call an Alexandrov-Fenchel type inequality because of its formal  resemblance to be an Alexandrov-Fenchel type inequality  for convex bodies.

\vskip 2mm
Following \cite{HardyLittlewoodPolya}, we say that two functions $f$ and $g$ are {\em effectively proportional}  if there are 
constants $a$ and $b$, not both zero, such that $af=bg$. Functions $f_1, \dots, f_m$ are effectively proportional if every pair $(f_i, f_j), 1\leq i, j\leq m$ is  effectively proportional.  A null function is effectively proportional to any function.
These notions will be used in the next theorems. 
\par
For a measure space 
$(X, \mu)$ and probability densities $p_i $ and
$q_i$, $1 \leq i \leq n$, we put
\begin{equation}\label{g0}
g_0(u)= \prod
_{i=1}^{n-m} \left[f_i\left(\frac{p_{i}}{q_{i}}\right)
q_{i}\right]^\frac{1}{n},
\end{equation} and for $j=0, \cdots, m-1$, 
\begin{equation}\label{gi}
g_{j+1}(u)=
\left[f_{n-j}\left(\frac{p_{n-j}}{q_{n-j}}\right)
q_{n-j}\right]^\frac{1}{n}.
\end{equation} For a vector $\vec{p}$, we denote by $\vec{p}^{\ n,k}$ the following vector $$\vec{p}^{\ n, k}=(p_1, \cdots, p_{n-m}, \underbrace{p_k, \cdots, p_k}_m), \ \ \ k>n-m.$$
\bt \label{inequality:mixed:f:divergence} Let $(X, \mu)$ be a measure space. For $1 \leq i \leq n$, let $P_i$ and
$Q_i$ be probability measures on $(X, \mu)$ with density functions $p_i$ and $q_i$ respectively  $\mu$-a.e. Let $f_i: (0,
\infty) \rightarrow \mathbb{R}^+$, $ 1 \leq i \leq n$, be convex
functions. Then, for $1 \leq m \leq n$, 
$$ \big[\dfpq\big]^m \leq \prod_{k =n-m+1}^ {n} D_{\vec{f}^{n, k}}\big(\vec{\bP}^{n, k},\vec{\bQ}^{n, k}\big).$$
  Equality holds if and only if one of the functions $g_0^\frac{1}{m} g_{i}$, $1 \leq i \leq m$,  is null or all are effectively proportional $\mu$-a.e. 
\par
\noindent If $m=n$,
\begin{eqnarray*} [\dfpq]^n
\leq \prod_{i=1}^n D_{f_i}(P_i, Q_i), 
\end{eqnarray*} 
with equality if and only if one of the functions  $f_{j}\left(\frac{p_{j}}{q_{j}}\right)
q_{j}$, $0 \leq j \leq n$, is null or all  are effectively proportional $\mu$-a.e.
\et

 \noindent {\bf Remarks.}  
(i) In particular, 
equality holds in Theorem
\ref{inequality:mixed:f:divergence} if all $(P_i, Q_i)$ coincide, and $f_i=\lambda_i f$ 
for some convex positive function $f$ and $\lambda_i\geq 0$, $i=1, 2, \cdots, n$.
\par
\noindent
(ii) 
 Theorem
\ref{inequality:mixed:f:divergence} still holds true if  the functions $f_i$ are
concave.

\vskip 2mm \noindent {\bf Proof.}  We let $g_0$ and $g_{j+1}$, $j=0, \cdots, m-1$
as in (\ref{g0}) and (\ref{gi}). 
By H\"{o}lder's inequality (see
\cite{HardyLittlewoodPolya}) 
\begin{eqnarray*}
[\dfpq]^m&=&\left(\int
_{X}g_0(u) g_1(u) \cdots g_{m}(u)\,d\mu\right)^m\\ &=& \bigg(\int
_{X}\prod_{j=0}^{m-1}\left[g_0(u) g_{j+1}(u)^m\right]^{\frac{1}{m}}\,d\mu\bigg)^m\\ &\leq& \prod _{j=0}^{m-1} \left(\int _{X} g_0(u)
g_{j+1}^m(u)\,d\mu\right)\\ &=&\prod_{k =n-m+1}^ {n} D_{\vec{f}^{n, k}}\big(\vec{\bP}^{n, k},\vec{\bQ}^{n, k}\big).
\end{eqnarray*}
  Equality holds in H\"{o}lder's inequality, if and only if 
one of the functions $g_0^\frac{1}{m} g_{i}$, $1 \leq i \leq m$,  is null or all are effectively proportional $\mu$-a.e. 
In particular, this is the case, if for all $i=1,  \cdots, n$, $(P_i, Q_i)=(P, Q)$ and $f_i=\lambda_i f$
for some convex function $f$ and $\lambda_i\geq 0$. 

\vskip 3mm We require some properties of $f$-divergences for our next result. 
Let $f:(0,\infty) \rightarrow \mathbb{R}^+$ be a convex function. By Jensen's inequality,
 \begin{equation}\label{Iso:type:1} D_f(P, Q)  =  \int_X
  f  \left(  \frac{p}{q}  \right)q\,d\mu\geq f  \left(  \int_X
  p\,d\mu  \right)=f(1),\end{equation} for all pairs of probability
measures $(P, Q)$ on $(X, \mu)$ with nonzero density functions $p$ and $q$ respectively $\mu$-a.e.  When $f$ is linear, equality holds trivially in (\ref {Iso:type:1}) . When $f$ is strictly convex, equality holds
true if and only if $p=q$ $\mu$-a.e.
If $f$ is a concave function,  Jensen's inequality implies   \begin{equation}\label{Iso:type:2}
D_f(P, Q)  =  \int_X   f  \left(  \frac{p}{q}  \right)q\,d\mu\leq
f  \left(  \int_X   p\,d\mu  \right)=f(1),\end{equation} for all pairs of
probability measures $(P, Q)$.  Again, when $f$ is linear, equality holds trivially. When $f$ is strictly concave,
equality holds true if and only if $p=q$ $\mu$-a.e.

\vskip 2mm For the mixed $f$-divergence with concave functions,
one has the following result.

\bt  \label{inequality:mixed:f:divergence2} 
Let $(X, \mu)$ be a measure space. For all $1 \leq i \leq n$, let
$P_i$ and $Q_i$ be probability measures on $X$ whose density functions $p_i$ and $q_i$ are nonzero $\mu$-a.e.
Let $f_i: (0, \infty) \rightarrow \mathbb{R}^+$, $ 1 \leq i \leq n$,
be concave functions. Then
\begin{eqnarray}\label{inequality:mixed:f:2} [\dfpq]^n
\leq \prod_{i=1}^n D_{f_i}(P_i, Q_i) \leq   \prod_{i=1}^nf_i(1).
\end{eqnarray} 

\noindent If in addition, all $f_i$ are strictly concave, equality holds if and only if there is a probability density $p$ such that for all $i=1, 2, \cdots n$,
 $$p_i=q_i=p, \  \  \mu-a.e.   $$ 
\et

\vskip 2mm 
\noindent {\bf Proof.} Theorem
\ref{inequality:mixed:f:divergence} and the remark after imply that for all concave
functions $f_i$,
\begin{eqnarray*}[\dfpq]^n \leq \prod_{i=1}^n D_{f_i}(P_i, Q_i)\leq  \prod_{i=1}^n f_i(1),
\end{eqnarray*} where the second inequality follows from
inequality (\ref{Iso:type:2}) and $f_i\geq 0$. 
\par
Suppose now that for all $i$, $p_i=q_i=p$, $\mu$-a.e., where $p$ is a fixed probability density. 
Then equality holds trivially in (\ref{inequality:mixed:f:2}).  Conversely, suppose that equality holds in (\ref{inequality:mixed:f:2}). Then, in particular, equality holds in Jensen's inequality which, as noted above, happens if and only if $p_i=q_i$ for all $i$. Thus,    $$\dfpq =\left(\prod_{i=1}^n [f_i(1)]^{1/n}\right) \int _{X} q_1^{1/n} \dots \,q_n^{1/n}d\mu.$$ Note also that if all $f_i: (0, \infty)\rightarrow \mathbb{R}^+$ are strictly concave, $f_{i}(1)\neq 0$ for all $1\leq i\leq n$.  Equality characterization in H\"older's inequality  implies that all $q_i$ are effectively proportional $\mu$-a.e. As all $q_i$ are probability measures,  they are all equal ($\mu$-a.e.) to a probability measure with density function (say) $p$. 
\vskip 2mm
\noindent
{\bf Remark.}
If $f_i(t)=a_i t +b_i$ are all linear and positive, then equality holds if and only if all $p_i, q_i$ are equal ($\mu$-a.e.) as convex combinations, i.e., if and only if for all $i, j$ $$
\frac{a_i}{a_i+b_i} p_i + \frac{b_i}{a_i+b_i} q_i = \frac{a_j}{a_j+b_j} p_j + \frac{b_j}{a_j+b_j} q_j,  \hskip 4mm \mu - \text{a.e.}
$$

\section{The $i$-th mixed $f$-divergence}
Let $(X, \mu)$ be a measure space. Throughout this section, we assume that the functions  $$f_1, f_2: (0, \infty)\rightarrow \{x\in \mathbb{R}: x>0\},$$ are  convex or concave, and that {$P_1, P_2, Q_1, Q_2$  are probability
measures on $X$ with density functions $p_1, p_2, q_1, q_2$ which are nonzero $\mu$-a.e.}  We also write $$\vf=(f_1, f_2), \ \ \vP=(P_1, P_2), \ \ \vQ=(Q_1, Q_2).$$

\begin{definition} Let $i\in \bbR$. The $i$-th mixed
$f$-divergence for $(\vf, \vP, \vQ)$,
denoted by $\dfpqi$,
is defined as \begin{equation} \dfpqi  =      \int_{X}    
\left[f_1  \left(\frac{p_{1}}{q_{1}}\right)   q_{1}\right]^\frac{i}{n}    
\left[f_2  \left(\frac{p_{2}}{q_{2}}\right)  
q_{2}\right]^\frac{n-i}{n}            d \mu.
\label{i:mixed:phi}\end{equation}
\end{definition} 
\vskip 2mm
\noindent
{\bf Remarks.} 
Note that  the  $i$-th mixed
$f$-divergence is  defined for any combination of convexity and concavity of $f_1$ and $f_2$,  namely,  both $f_1$ and $f_2$ concave, or both $f_1$ and $ f_2$ convex, or one is  convex
the other  is concave.  
\par
 
It is easily checked that $$\dfpqi=D_{(f_2, f_1)}\big((P_2, Q_2), (P_1,
Q_1); n-i\big).$$
If $0\leq i\leq n$ is an integer, then the triple  $(f_1,
P_1, Q_1)$ appears $i$-times while the triple $(f_2, P_2, Q_2)$ appears
$(n-i)$ times in $\dfpqi$.  
Note that if $i=0$, then  
 $\dfpqi=D_{f_2}(P_2,
Q_2),$ and if $i=n$ then  $\dfpqi=D_{f_1}(P_1, Q_1).$ 
\par
Another special case is when
$P_2=Q_2=\mu$ almost everywhere and $\mu$ is also a probability measure. Then such an $i$-th mixed $f$-divergence, denoted by
$D\big((f_1, P_1, Q_1), i; f_2\big)$, has the  form
\begin{equation*}
D\big((f_1, P_1, Q_1), i; f_2\big)=[f_2(1)]^{1-i/n}
\int_{X} \left[f_1\left(\frac{p_{1}}{q_{1}}\right)
q_{1}\right]^\frac{i}{n} d \mu.
\end{equation*}
 \vskip 2mm
\noindent
{\bf Examples and Applications.}
\par
\noindent (i) For $f(t)=|t-1|$, we get the {\em $i$-th mixed total variation} $$D_{TV}\big(\vP, \vQ; i\big) =\int_X |p_1-	q_1|^{\frac{i}{n}} |p_2-	q_2|^{\frac{n-i}{n}}\,d\mu.$$  
\par
\noindent
(ii) For $f_1(t)= f_2(t)= [t \ln t]_+$, we get the (modified)  {\em $i$-th mixed relative entropy} or {\em
 $i$-th mixed Kullback Leibler divergence}
 \begin{eqnarray*}
D_{KL}\big(\vP, \vQ; i\big)  = \int_{X} 
\left[ p_1 \ln\left(\frac{p_{1}}{q_{1}}\right) 
\right]_+^\frac{i}{n}  \left[ p_2  \ln\left(\frac{p_{2}}{q_{2}}\right) 
\right]_+^\frac{n-i}{n}  d\mu.
\end{eqnarray*}
\par
\noindent
(iii) For the convex or concave functions $f_{\alpha_j}(t) = t^{\alpha_j}$, $j=1,2$, we get the {\em $i$-th mixed Hellinger integrals}
\begin{eqnarray*}
D_{(f_{\alpha_1}, f_{\alpha_2})}\big( \vP, \vQ; i\big) = \int_{X} 
\left( p_1  ^{\alpha_1} q_1 ^{1-\alpha_1} \right) ^\frac{i}{n}   \left( p_2  ^{\alpha_2} q_2 ^{1-\alpha_2} \right) ^\frac{n-i}{n} 
 d\mu.
\end{eqnarray*}
In particular, for $\alpha_j= \alpha$, for  $j=1,2$, 
\begin{eqnarray*}
D_{(f_{\alpha}, f_{\alpha})}\big( \vP, \vQ; i\big)  = \int_{X} 
\left( p_1  ^{\alpha} q_1 ^{1-\alpha} \right) ^\frac{i}{n}   \left( p_2  ^{\alpha} q_2 ^{1-\alpha} \right) ^\frac{n-i}{n} 
 d\mu.
\end{eqnarray*}
This integral  can be used to define the {\em $i$-th mixed $\alpha$-R\'enyi divergence }
 \begin{eqnarray*}
D_{\alpha}\big(\vP, \vQ; i \big) = \frac{1}{\alpha -1}  \ln \left[D_{(f_{\alpha}, f_{\alpha})}\big( \vP, \vQ; i\big)\right] .
\end{eqnarray*}
The case $\alpha_i=\frac{1}{2}$ for all $i$ gives 
\begin{eqnarray*}
D_{(\sqrt{t}, \sqrt{t})}\big( \vP, \vQ; i\big)  = \int_{X} 
\left( p_1  q_1  \right) ^\frac{i}{2n}   \left( p_2   q_2 \right) ^\frac{n-i}{2n} 
 d\mu,
\end{eqnarray*}
the {\em $i$-th mixed Bhattacharyya  coefficient} or 
{\em $i$-th mixed Bhattacharyya distance} of $p_i$ and $q_i$.
\par
\noindent
(iv) Important  applications are again in the theory of convex bodies.  As in section 2,
 let $K_1$ and $K_2$ be convex bodies with positive curvature function. For $l=1,2$, let
\begin{equation*}
\,d P_{K_l}= \frac{1}{h_{K_l}^n}\,d\sigma  \ \ \ \text{and}   \ \ \   \,d Q_{K_l}={f_{K_l} h_{K_l}}\,d\sigma. 
\end{equation*} Let $f_l: (0, \infty) \rightarrow \mathbb{R}$, $ l=1, 2$, be positive convex functions.
Then, we define the {\em $i$-th mixed $f$-divergence} for convex bodies $K_1$ and $ K_2$ 
by
\begin{eqnarray*}
 D_{\vf}\big( (P_{K_1}, P_{K_2}), (Q_{K_1}, Q_{K_2}); i\big)= \int_{S^{n-1}} \bigg[f_1 \bigg(\frac{1}{ f_{K_1}h_{K_1}^{n+1}} \bigg)   {f_{K_1} h_{K_1}}   \bigg]^\frac{i}{n}    \bigg [f_2 \bigg(\frac{1}{f_{K_2}h_{K_2}^{n+1}} \bigg)  {f_{K_2} h_{K_2}}  \bigg]^\frac{n-i}{n}   \, d\sigma.
\end{eqnarray*}
 These are the general $i$-th mixed affine surface areas introduced in \cite{Ye2012}.

\vskip 2mm  The following result holds for all possible combinations of convexity and concavity 
of $f_1$ and $ f_2$.
\vskip 2mm
\begin{proposition} \label{Monotone:1} Let $\vf, \vP, \vQ$ be as above. 
If $j\leq i\leq k$ or $k\leq i\leq j$, then
\begin{eqnarray*} \dfpqi\leq  \bigg[D_{\vf}\big(\vP, \vQ; j\big)\bigg]^{\frac{k-i}{k-j}}\times \bigg[D_{\vf}\big(\vP, \vQ; k\big)\bigg]^{\frac{i-j}{k-j}}.\end{eqnarray*}
Equality holds trivially if $i=k$ or $i=j$. Otherwise, 
equality holds if and only if one of the functions $f_i\left(\frac{p_i}{q_i}\right) q_i$, $i=1,2$,  is null,  or $f_1\left(\frac{p_1}{q_1}\right) q_1$ and $f_2\left(\frac{p_2}{q_2}\right) q_2$ are effectively proportional $\mu$-a.e.
In particular, this holds if $(P_1, Q_1)=(P_2, Q_2)$ and $f_1=\l f_2$ for
some $\l>0$.
\end{proposition}

\vskip 2mm \noindent {\bf Proof.} By formula (\ref{i:mixed:phi}),
one has  \begin{eqnarray*} \dfpqi    &    =    &    \int_{X}
\left[f_1\left(\frac{p_{1}}{q_{1}}\right) q_{1}\right]^\frac{i}{n}
\left[f_2\left(\frac{p_{2}}{q_{2}}\right)
q_{2}\right]^\frac{n-i}{n} d \mu\\     &    =    &    \int
_{X}    \left\{  \left[   f_1\left(  \frac{p_{1}}{q_{1}}  \right)
q_{1}  \right]^\frac{j}{n}     \left[  f_2  \left(  \frac{p_{2}}{q_{2}}  \right)
q_{2}  \right]^\frac{n-j}{n}  \right\}^{\frac{k-i}{k-j}}\\ && \times  \left\{  \left[  f_1\left(  \frac{p_{1}}{q_{1}}  \right)
q_{1}  \right]^\frac{k}{n}     \left[  f_2\left(  \frac{p_{2}}{q_{2}}  \right)
q_{2}  \right]^\frac{n-k}{n}  \right\}^{\frac{i-j}{k-j}}         d\mu \\     &    \leq    &     \bigg[D_{\vf}\big(\vP, \vQ; j\big)\bigg]^{\frac{k-i}{k-j}}\times \bigg[D_{\vf}\big(\vP, \vQ; k\big)\bigg]^{\frac{i-j}{k-j}},
\end{eqnarray*} 
where the last inequality follows from H\"{o}lder's  inequality and
formula (\ref{i:mixed:phi}). 
The equality characterization follows from the one in H\"{o}lder inequality.
In particular, if $(P_1, Q_1)=(P_2, Q_2)$,
and $f_1=\l f_2$ for some $\l>0$,  equality holds. 
\vskip 3mm
\bc \label{KOR}Let $f_1$ and $ f_2$ be positive, concave functions on $(0, \infty)$. Then
for all $\vP, \vQ$ and  for all
$0\leq i\leq n$,
\begin{equation*} \big[\dfpqi\big]^n\leq [f_1(1)]^i
[f_2(1)]^{n-i}.\end{equation*}
If in addition, $f_1$ and $f_2$ are strictly concave, equality holds iff 
$p_1=p_2=q_1=q_2$ $\mu$-a.e.
\ec
\vskip 2mm 
\noindent {\bf Proof.} Let $j=0$ and  $k=n$ in Proposition
\ref{Monotone:1}.  Then for all $0\leq i\leq n$,
\begin{eqnarray*}\label{i:mixed:phi:1}
\big[\dfpqi\big]^n \leq
[D_{f_1}(P_1, Q_1)]^{i}[D_{f_2}(P_2, Q_2)]^{n-i} \leq  [f_1(1)]^i [f_2(1)]^{n-i},\nonumber\end{eqnarray*} where the
last inequality follows from inequality (\ref{Iso:type:2}).
\par
To have equality, the above inequalities should be equalities.  Proposition \ref{Monotone:1} implies that  $f_1\left(\frac{p_1}{q_1}\right) q_1$ and $f_2\left(\frac{p_2}{q_2}\right) q_2$ are effectively proportional $\mu$-a.e. As both  $f_1$ and $f_2$ are strictly concave,  Jensen's inequality requires that $p_1=q_1$ and $p_2=q_2$ $\mu$-a.e. Therefore, equality holds if and only if $f_1(1)q_1$ and $f_2(1)q_2$ are effectively proportional $\mu$-a.e. As both  $f_1(1)$ and $f_2(1)$  are not zero,  equality holds iff $p_1=p_2=q_1=q_2$ $\mu$-a.e.   

\vskip 2mm
\noindent
{\bf Remark.}
If $f_1(t) =a_1t+b_1$ and $f_2(t)=a_2t+b_2$ are both linear,  equality holds in Corollary \ref{KOR}  if and only if
$p_i, q_i$, $i=1,2$, are equal as convex combinations,  i.e.,  $$
\frac{a_1}{a_1+b_1} p_1 + \frac{b_1}{a_1+b_1} q_1 = \frac{a_2}{a_2+b_2} p_2 + \frac{b_2}{a_2+b_2} q_2,  \hskip 4mm \mu - \text{a.e.}
$$
\vskip 2mm

This proof can be used to establish the following result for
$D\big((f_1, P_1, Q_1), i; f_2\big)$.

\bc Let $(X, \mu)$ be a probability space. Let $f_1$ be a positive concave
function on $(0, \infty)$. Then for all $P_1, Q_1$, for
all (concave or convex) positive functions $f_2$, and for all $0\leq
i\leq n$,
\begin{equation*} \big[D\big((f_1, P_1, Q_1), i;
f_2\big)\big]^n\leq [f_1(1)]^i [f_2(1)]^{n-i}.\end{equation*} 
If $f_1$ is strictly concave, equality holds  if and only if $P_1=Q_1=\mu$.
When $f_1(t) =at+b$ is linear, equality holds if and only if ${ap_1+bq_1}={a+b}$ $\mu$-a.e.\ec

\vskip 2mm

\bc \label{KOR1} Let $f_1$ be a positive convex function and $f_2$ be a positive concave function on $(0,
\infty)$. Then, for all $\vP, \vQ$, and  for all $k\geq n$,
\begin{equation*} \big[D_{\vf}\big(\vP, \vQ;  k\big)\big]^n\geq [f_1(1)]^k
[f_2(1)]^{n-k}.\end{equation*} 
If in addition, $f_1$ is strictly convex and $f_2$ is strictly concave, equality holds
if and only if $p_1=p_2=q_1=q_2$ $\mu$-a.e.
\ec 
\vskip 2mm 
\noindent 
{\bf Proof.}  On the right hand side of Proposition \ref{Monotone:1}, let $i=n$ and $ j=0$. Let  $k\geq n$. Then
\begin{eqnarray*}
\big[D_{\vf}\big(\vP, \vQ;  k\big)\big]^n \geq  [D_{f_1}(P_1, Q_1 )]^{k}[D_{f_2}(P_2, Q_2
)]^{n-k} \geq  [f_1(1)]^k [f_2(1)]^{n-k}.\end{eqnarray*} Here, the
last inequality follows from inequalities (\ref{Iso:type:1}), (\ref{Iso:type:2}) and $k\geq n$. 
To have equality, the above inequalities should be equalities.  Proposition \ref{Monotone:1} implies  that  $f_1\left(\frac{p_1}{q_1}\right) q_1$ and $f_2\left(\frac{p_2}{q_2}\right) q_2$ are effectively proportional $\mu$-a.e. As $f_1$ is strictly convex and $f_2$ is strictly concave,  Jensen's inequality implies that $p_1=q_1$ and $p_2=q_2$ $\mu$-a.e. Therefore, as 
both $f_1(1)$ and $f_2(1)$ are not zero, equality holds if and only if $p_1=p_2=q_1=q_2$ $\mu$-a.e.  
 
\vskip 2mm
\noindent
{\bf Remark.}
If $f_1(t) =a_1t+b_1$ and $f_2(t)=a_2t+b_2$ are both linear,  equality holds in Corollary \ref{KOR1}  if and only if
$p_i, q_i$, $i=1,2$, are equal $\mu$-a.e. as convex combinations, i.e.,  $$
\frac{a_1}{a_1+b_1} p_1 + \frac{b_1}{a_1+b_1} q_1 = \frac{a_2}{a_2+b_2} p_2 + \frac{b_2}{a_2+b_2} q_2,  \hskip 4mm \mu - \text{a.e.}
$$

This proof can be used to establish the following result for
$D\big((f_1, P_1, Q_1), k; f_2\big)$. 
\vskip 3mm
\bc Let $(X, \mu)$ be a probability space. Let $f_1$ be a positive convex
function on $(0, \infty)$. Then for all $P_1, Q_1$,  for
all (positive concave or convex) functions $f_2$, and for all $k\geq n$,
\begin{equation*} \big[D\big((f_1, P_1, Q_1), k;
f_2\big)\big]^n\geq [f_1(1)]^k [f_2(1)]^{n-k}. \end{equation*}If $f_1$ is strictly convex, equality holds if and only if $P_1=Q_1=\mu$.
When $f_1(t) =at+b$ is linear, equality holds if and only if ${ap_1+bq_1}={a+b}$ $\mu$-a.e.

\ec
\vskip 3mm
\bc 
Let $f_1$ be a positive concave function and $f_2$ be a positive convex function on $(0, \infty)$. Then for
all $\vP, \vQ$,  and  for all $k\leq 0$,
\begin{equation*} \big[D_{\vf}(\vP, \vQ; k)\big]^n\geq [f_1(1)]^k [f_2(1)]^{n-k}.\end{equation*} 
If in addition, $f_1$ is strictly concave and $f_2$ is strictly convex, equality holds iff 
$p_1=p_2=q_1=q_2$ $\mu$-a.e.
\ec
\vskip 2mm \noindent {\bf Proof.}  Let $i=0$ and $j=n$ in Proposition
\ref{Monotone:1}. Then
\begin{eqnarray*}
\big[D_{\vf}\big(\vP, \vQ; k)\big]^n  &\geq&  [D_{f_1}(P_1, Q_1 )]^{k}[D_{f_2}(P_2, Q_2 )]^{n-k} \geq  [f_1(1)]^k [f_2(1)]^{n-k}.\end{eqnarray*} Here, the last
inequality follows from inequalities (\ref{Iso:type:1}), (\ref{Iso:type:2}), and $k\leq 0$.

To have equality, the above inequalities should be equalities. Proposition \ref{Monotone:1} implies that  $f_1\left(\frac{p_1}{q_1}\right) q_1$ and $f_2\left(\frac{p_2}{q_2}\right) q_2$ are effectively proportional $\mu$-a.e. As $f_1$ is strictly concave and $f_2$ is strictly convex, Jensen's inequality requires that $p_1=q_1$ and $p_2=q_2$. Therefore, equality holds if and only if $f_1(1)q_1$ and $f_2(1)q_2$ are effectively proportional $\mu$-a.e.  As both $f_1(1)$ and $f_2(1)$ are not zero,  equality holds if and only if $p_1=p_2=q_1=q_2$ $\mu$-a.e.  
 
 \vskip 3mm 
This proof can be used to establish the following result for $D\big((f_1, P_1, Q_1), k; f_2\big)$.
 \vskip 2mm
\bc 
Let $f_1$ be a concave function on $(0, \infty)$. Then for
all $P_1, Q_1$, for all (concave or convex) functions
$f_2$, and for all $k\leq 0$,
\begin{equation*} \big[D\big((f_1, P_1, Q_1), k;  f_2\big)\big]^n\geq [f_1(1)]^k [f_2(1)]^{n-k}.\end{equation*}
If $f_1$ is strictly concave, equality holds if and only if $P_1=Q_1=\mu$.
When $f_1(t) =at+b$ is linear, equality holds if and only if ${ap_1+bq_1}={a+b}$ $\mu$-a.e.

\ec

 \vskip 2mm \noindent Elisabeth Werner, {\small \tt elisabeth.werner@case.edu}\\
 {\small Department of Mathematics \ \ \ \ \ \ \ \ \ \ \ \ \ \ \ \ \ \ \ Universit\'{e} de Lille 1}\\
 {\small Case Western Reserve University \ \ \ \ \ \ \ \ \ \ \ \ \ UFR de Math\'{e}matique }\\
 {\small Cleveland, Ohio 44106, U. S. A. \ \ \ \ \ \ \ \ \ \ \ \ \ \ \ 59655 Villeneuve d'Ascq, France}\\
  \\
 \and Deping Ye,  {\small \tt deping.ye@mun.ca}\\
 {\small Department of Mathematics and Statistics}\\
 {\small  Memorial University of Newfoundland}\\
 {\small St. John's, Newfoundland, Canada A1C 5S7}\\
 
\end{document}